\documentclass[12pt]{amsart}
\usepackage{amsmath}%,amsthm}%,amssymb
\usepackage{amscd}
\usepackage[psamsfonts]{amssymb}
\usepackage[mathscr,psamsfonts]{eucal}
\DeclareMathAlphabet{\eurm}{U}{eur}{m}{n}
\DeclareMathAlphabet{\cyrm}{U}{UWCyr}{m}{n}
\newcommand{\M}[1]{{\mathscr{#1}}}

\newcommand{\f}[1]{{\boldsymbol{#1}}}

\newcommand{\K}[1]{{\cyrm{#1}}}
\DeclareMathOperator{\id}{{id}}
\DeclareMathOperator{\Ob}{{Ob}}

\DeclareMathOperator{\Lin}{{Lin}}

\DeclareMathOperator{\Cla}{Cla}
\DeclareMathOperator{\pr}{pr}
\DeclareMathOperator{\byd}{\,{\raisebox{.1ex}{$\eurm :$}{\eurm =}}\,}
\newcommand{\sig}{\sigma}\newcommand{\del}{\delta}
\newcommand{\alp}{\alpha}
\newcommand{\bet}{\beta}\newcommand{\eps}{\epsilon}
\newcommand{\lam}{\lambda}

\newcommand{\Lam}{\Lambda}
\newcommand{\Gam}{\Gamma}
\newcommand{\gam}{\gamma}
\newcommand{\bEq}{\begin{eqnarray}}
\newcommand{\eEq}{\end{eqnarray}}
\newcommand{\beq}{\begin{eqnarray*}}
\newcommand{\eeq}{\end{eqnarray*}}
\newcommand{\dt}[1]{{\dot{#1}}}
\newcommand{\uten}[1]{\underset{#1}{\otimes}}
\newcommand{\ucar}[1]{\underset{#1}{\times}}
\newcommand{\sep}[1]{\qquad\text{#1}\qquad}
\newcommand{\ba}[1]{{{\bar{#1}}}}
\newcommand{\col}[3]{_{#1}{}^{#2}{}_{#3}}
\newcommand{\der}{\partial}
\newcommand{\nab}{\nabla}

\newcommand{\com}{\circ}
\newcommand{\car}{\times}
\newcommand{\ten}{\otimes}
\newcommand{\wed}{\wedge}
\newcommand{\Rn}{{\mathbb R}}

\newcommand{\tfr}[2]{\tfrac{#1}{#2}\,}
\newcommand{\END}{{\,\text{\qedsymbol}}}

\begin{document}

\newcounter{theorem}

\newtheorem{definition}[theorem]{Definition}
\newtheorem{lemma}[theorem]{Lemma}
\newtheorem{proposition}[theorem]{Proposition}
\newtheorem{theorem}[theorem]{Theorem}
\newtheorem{corollary}[theorem]{Corollary}
\newtheorem{remark}[theorem]{Remark}
\newtheorem{example}[theorem]{Example}
\newtheorem{Note}[theorem]{Note}
\newcounter{assump}
\newtheorem{Assumption}{\indent Assumption}[assump]
\renewcommand{\thetheorem}{\thesection.\arabic{theorem}}

\newcommand{\bCr}{\begin{corollary}}
\newcommand{\eCr}{\end{corollary}}
\newcommand{\bDf}{\begin{definition}\em}
\newcommand{\eDf}{\end{definition}}
\newcommand{\bLm}{\begin{lemma}}
\newcommand{\eLm}{\end{lemma}}
\newcommand{\bPr}{\begin{proposition}}
\newcommand{\ePr}{\end{proposition}}
\newcommand{\bRm}{\begin{remark}\em}
\newcommand{\eRm}{\end{remark}}
\newcommand{\bEx}{\begin{example}\em}
\newcommand{\eEx}{\end{example}}
\newcommand{\bTh}{\begin{theorem}}
\newcommand{\eTh}{\end{theorem}}
\newcommand{\bNt}{\begin{Note}\em}
\newcommand{\eNt}{\end{Note}}
\newcommand{\bPf}{\begin{proof}[\noindent\indent{\sc Proof}]}
\newcommand{\ePf}{\end{proof}}
%--------------------------------------------------------------------%

\title[Natural connections]
        {Natural connections given by general
        linear\\ and classical connections}

\author{Josef Jany\v ska}



\address{
\newline
{Department of Mathematics, Masaryk University
\newline
Jan\'a\v ckovo n\'am. 2a, 662 95 Brno, Czech Republic}
\newline
E-mail: {\tt janyska@math.muni.cz}
}
\thanks{This paper has been supported
by the Grant agency of the Czech Republic under the Project
number GA 201/02/0225.}

%--------------------------------------------------------------------%
\begin{abstract}
We assume a vector bundle $p:\f E\to \f M$ with a general linear
connection $K$ and a classical linear connection $\Lam$ on $\f M$.
We prove that all classical linear connections on the total space
$\f E$ naturally given by $(\Lam, K)$ form a 15-parameter family.
Further we prove that all connections on $J^1\f E$
naturally given by $(\Lam, K)$ form a 14-parameter family.
Both families of connections are described geometrically.
\end{abstract}
\maketitle
%--------------------------------------------------------------------%

\section*{Introduction}
\setcounter{equation}{0}
\setcounter{theorem}{0}
%--------------------------------------------------------------------%

The reduction theorems for classical (linear) connections
on manifolds are very powerful tools to classify natural
(invariant) operators (of any finite order) from the
bundle of classical connections and
some tensor bundle to another tensor bundle. By reduction theorems
all such operators reduce to operators defined on tensor bundles
only (operators of curvature tensor fields of classical connections
in question, given tensor fields and covariant
differentials of the curvature
and given tensor fields). For the
general theory of natural bundles  and operators see
\cite{KolMicSlo93, KruJan90, Nij72, Ter78} and
for the proof of reduction theorems for classical connections and
their applications see
\cite{KolMicSlo93, Sch54}.

The reduction idea for  general  linear connections on vector bundles
is also possible.  It is a gauge version of reduction theorems
for the gauge group
$G=GL(n,\Rn)$. For operators of order one the first reduction theorem
for general linear connection is in fact the special case
of the Utiyama's theorem,
see \cite{Eck81, KolMicSlo93}. In \cite{Jan04a} the reduction
theorems for general linear connections were proved
for operators of any finite orders with values in a gauge-natural
bundle of order (1,0). In this gauge situation auxiliary
classical connections on base manifolds have to be used.

Sometimes we need to study natural operators which have values
in a natural or  a gauge-natural bundle of higher order.
In this case we can use
higher order valued reduction theorems,
for the case of classical connections
see \cite{Jan04b} and for the case of general linear connections
see \cite{ Jan04c}. Typical
operators of this type are natural tensor fields on the tangent
or the cotangent bundle of a manifold with a classical connection
or natural tensor fields
on the total space of a vector bundle with a general linear connection
and a classical connection on the base manifold.

In this paper we use the higher order valued reduction theorems
for general linear connections, \cite{Jan04c}, to clasify
all classical linear connections on the total space
of a vector bundle $p:\f E\to \f M$
and connections on the first jet prolongation $J^1\f E$ naturally
induced by a general linear connection $K$ on $\f E$ and a
classical connection $\Lam$ on $\f M$. We shall prove that there is
a 15-parameter family of classical linear connections
$\widetilde{D}(\Lam, K)$ on
the total space $\f E$ and a 14-parameter family of connections
$\widetilde{\Gam}(\Lam,K)$ on $J^1\f E$. We present also the geometrical
description of both families and we show that
$\widetilde{\Gam}(\Lam,K)$
can be obtained from $\widetilde{D}(\Lam, K)$ by using
the operator $\chi$ described
in \cite{JanMod96}.  The result obtained for classical linear
connections on $\f E$ coincides with the classification of
\cite{GanKol91} (see also \cite{KolMicSlo93}).

All manifolds and maps are assumed to be smooth.
The sheaf of (local) sections of a fibered  manifold $p:\f Y\to \f X$
is denoted by $C^\infty(\f Y)$, $C^\infty(\f Y,\Rn)$ denotes the sheaf
of (local) functions.

%--------------------------------------------------------------------%
\section{Higher order valued reduction theorems for general linear
        connections}
\setcounter{equation}{0}
\setcounter{theorem}{0}
%--------------------------------------------------------------------%

In what follows let $G=GL(n,\Rn)$ be the group of linear automorphisms
of $\Rn^n$ with coordinates $(a^i_j)$, $i,j=1,\dots,n$.
Let $\M M_m$  be the category of $m$-dimensional
$C^\infty$-manifolds  and  smooth  embeddings. Let $\M{VB}_{m,n}$
be the category
of vector bundles with $m$-dimensional bases,
$n$-dimensional fibers and local fiber linear diffeomorphisms
and
let $\M P\M B_m(G)$  be
the category of smooth  principal  $G$-bundles  with  $m$-dimensional
bases  and
smooth $G$-bundle maps $(\varphi,f)$, where the map $f \in$
Mor$\M M_m$.
Then any vector bundle
$
(p: \f E\to \f M) \in\Ob\M{VB}_{m,n}
$
can be considered as a zero
order $G$-gauge-natural bundle given by the $G$-gauge-natural
bundle functor
of associated bundles
$\M{PB}_{m}(G)\to \M{VB}_{m,n}$.

Local linear fiber coordinates on $\f E$ will be denoted
by $(x^A)=(x^\lam, y^i)$,
$A=1,\dots,m+n$, $\lam=1,\dots,m$, $i=1,\dots,n$.
The induced fiber coordinates on
$T\f E$ or $T^*\f E$
will be denoted by
$(x^\lam, y^i,\dot x^\lam,\dot y^i)$
or $(x^\lam, y^i, \dot x_\lam,\dot y_i)$ and the induced
local bases of sections of
$T\f E$ or
$T^*\f E$
will be denoted by
$(\der_\lam,\der_i)$ or $(d^\lam, d^i)$, respectively.

\smallskip

We define a {\em  general linear connection\/}
on $\f E$ to be a linear section
$
K : \f E \to J^1\f E\,.
$
Considering the contact morphism $J^1\f E\to
T^*\f M \ten T\f E$ over the identity of $T\f M$, a
linear connection can be regarded as a $T\f E$-valued
1-form
$
K : \f E \to T^*\f M \ten T\f E
$
projecting onto the identity of $T\f M$.

The coordinate expression of a linear connection
$K$
is of the type
\bEq
K = d^\lam \ten \big(\der_\lam
+ K\col ji\lam \, y^j \, \der_i \big) \,,
\sep{with}  K\col ji \lam \in C^\infty(\f M,\Rn)\,.
\eEq

Linear connections can be regarded as sections of
a (1,1)-order $G$-gauge-natural bundle
$\Lin\f E \to \f M$, \cite{Eck81, Jan04a, KolMicSlo93}.
The standard fiber of the functor $\Lin$ will be denoted by
${R}=\Rn^{n*}\ten \Rn^n\ten\Rn^{m*}$, elements of $R$ will be
said to be {\em formal linear connections},
the induced coordinates
on ${R}$ will be said to be {\em formal symbols} of formal linear
connections and will be denoted by $(K\col ji\lam)$.

The curvature tensor field
$
R[{K}] :\f M\to\f E^* \ten \f E \ten \bigwedge^2 T^*\f M\,
$
given by $K$ is a natural operator
$
R[{K}]:C^\infty(\Lin\f E) \to
        C^\infty(\f E^*\ten\f E\ten\Lam^2T^*\f M)
$
which is of order one.
Its coordinate
expression is
\begin{align*}
R[{K}] & = R[K]\col ji{\lam\mu}\, d^j\ten  \der_i \ten d^\lam \wed d^\mu
\\
& = -2(\der_\lam K \col ji\mu + K\col jp\lam K\col pi\mu)\, d^j\ten
        \der_i \ten d^\lam \wed d^\mu \,.
\end{align*}

\smallskip

We define a {\em classical connection\/} on $\f M$ to be
a linear connection on the tangent vector bundle
$p_{\f M}:T\f M\to \f M$ with the coordinate
expression
\bEq
\Lam = d^\lam \ten \big(\der_\lam
+ \Lam\col\nu\mu\lam \, \dot x^\nu \, \dt\der_\mu \big) \,,
\quad  \Lam\col \mu \lam \nu\in C^\infty(\f M,\Rn)\,.
\eEq

Classical connections can be regarded as sections of a 2nd
order natural bundle
$\Cla\f M \to \f M$, \cite{KolMicSlo93}.
The standard fiber of the functor $\Cla$ will be denoted by
${Q}=\Rn^{m*}\ten \Rn^m\ten\Rn^{m*}$, elements of $Q$ will be said
to be {\em formal classical connections},
the induced
coordinates on ${Q}$ will be said to be {\em formal Christoffel
symbols} of formal classical connections and will be denoted by
$(\Lam\col \mu \lam \nu)$. The natural subbundle of symmetric
classical connections will be denoted by $\Cla_\tau \f M\to \f M$
with local fiber coordinates $(x^\lam, \Lam\col\mu\lam\nu)$,
$\Lam\col\mu\lam\nu= \Lam\col\nu\lam\mu$.

\smallskip

The curvature tensor field of a classical connection
is a natural operator
\beq
R[\Lam]:C^\infty(\Cla\f M) \to C^\infty(T^*\f M\ten T\f M\ten
        \bigwedge^2T^*\f M)
\eeq
which is of order one.

\smallskip

Let us denote by
$
\f E^{p,r}_{q,s}\byd \ten^p \f E
\ten\ten^q \f E^*\ten\ten^rT\f M \ten \ten^s T^*\f M
$
the tensor product over $\f M$ and recall that
$
\f E^{p,r}_{q,s}
$
is a vector bundle which is a $G$-gauge-natural bundle of order
$(1,0)$.

A classical connection $\Lam$ on $\f M$ and a linear
connection ${K}$ on $\f E$ induce the linear tensor product connection
$K^p_q\ten {\Lam}^r_s\byd \ten^p K\ten\ten^q K^*\ten
\ten^r{\Lam} \ten \ten^s{\Lam}^*$ on
$\f E^{p,r}_{q,s}$
\beq
K^p_q\ten {\Lam}^r_s: \f E^{p,r}_{q,s} \to T^*\f M\uten{\f M}T
        \f E^{p,r}_{q,s}
\eeq
which can be considered
as a linear section
\beq
K^p_q\ten{\Lam}^r_s:\f E^{p,r}_{q,s}\to J^1\f E^{p,r}_{q,s}\,.
\eeq

Then we define, \cite{Jan03}, the
{\em covariant differential of a section $\Phi:\f M\to \f E^{p,r}_{q,s}$
with respect to
the pair of connections $(\Lam,K)$} as a section of
$\f E^{p,r}_{q,s}\ten T^*\f M$ given by
\beq
\nab^{(\Lam,K)} \Phi = j^1 \Phi - (K^p_q\ten{\Lam}^r_s )\com \Phi\,.
\quad
\eeq

In what follows we set $\nab = \nab^{(\Lam,K)}$ and
$\phi^{i_1\dots i_p\lam_1\dots\lam_r}
        _{j_1 \dots j_q\mu_1\dots\mu_s;\nu}=
\nab_{\nu}\phi^{i_1\dots i_p\lam_1\dots\lam_r}
        _{j_1 \dots j_q\mu_1\dots\mu_s}$.
We shall denote by $\nab^i$ the $i$-th iterated covariant differential
and we shall put $\nab^{(k,r)}\byd(\nab^k,\dots,\nab^r)$, $r\ge k$,
$\nab^{(r)}\byd \nab^{(0,r)}$.

\smallskip

Let us denote by $C_{C,i}\f M$ the $i$-th {\em curvature bundle
of symmetric classical connections\/} given as
the image of the natural operator
$$
\nab^i R[\Lam]:J^{i+1}\Cla_{\tau}\f M\to
        T^*\f M\uten{\f M} T\f M\uten{\f M}\bigwedge^2T^*\f M
        \uten{\f M}\ten^i T^*\f M\,.
$$
We shall put $C^{(r)}_C\f M\byd C_{C,0}\f M\ucar{\f M}\dots
\ucar{\f M} C_{C,r}\f M$ and $\pr^r_k:C^{(r)}_C\f M\to
C^{(k)}_C\f M$, $r\ge k$, the canonical projection.

Similarly let us denote by $C_{L,i}\f E$ the $i$-th
{\em curvature bundle of general linear connections\/} given as
the image of the natural operator
$$
\nab^i R[K]:J^{i-1}\Cla_{\tau}\f M\ucar{\f M}J^{i+1}\Lin \f E\to
\f E^*\uten{\f M} \f E\uten{\f M}\bigwedge^2T^*\f M\uten{\f M}
        \ten^i T^*\f M\,.
$$
We shall put $C^{(r)}_L\f E\byd C_{L,0}\f E\ucar{\f M}\dots
\ucar{\f M} C_{L,r}\f E$ and $\pr^r_k:C^{(r)}_L\f M\to
C^{(k)}_L\f M$, $r\ge k$, the canonical projection.

We denote by
$
J^{k-2}\Cla_{\tau} \f M\ucar{\f M} J^{k-1} \Lin \f E\ucar{\f M}
C^{(k-2,s-1)}_C\f M\ucar{\f M} C^{(k-1,r-1)}_L\f E
$, $s\ge r-2$,
the pullback of
$$
\pr^{s-1}_{k-3}\car\pr^{r-1}_{k-2}:C^{(s-1)}_C\f M\ucar{\f M}
        C^{(r-1)}_L\f E\to C^{(k-3)}_C\f M\ucar{\f M}
        C^{(k-2)}_L\f E
$$
with respect to the surjective submersion, \cite{Jan04a},
$$
(\nab^{(k-3)}R[\Lam],\nab^{(k-2)}R[K]):J^{k-2}\Cla_{\tau}\f M
        \ucar{\f M} J^{k-1}\Lin\f E\to C^{(k-3)}_C\ucar{\f M}
        C^{(k-2)}_L\f E\,.
$$
Then the first $k$-th order valued reduction theorem can be
formulated as follows, for the poof see \cite{Jan04c},

\bTh\label{Th1.1}
Let $s\ge r-2$, $r+1, s+2\ge k \ge 1$.
Let $F$ be a $G$-gauge-natural bundle of order $k$.
All natural differential operators
$$
f:C^\infty (\Cla_{\tau}\f M\ucar{\f M} \Lin \f E) \to  C^\infty(F\f E)
$$
which are of order $s$ with respect to symmetric
classical connections and of order $r$ with respect
to general linear connections are of the
form
$$
f(j^s\Lam, j^r K) = g(j^{k-2}\Lam, j^{k-1} K,
\nab^{(k-2,s-1)} R[\Lam], \nab^{(k-1,r-1)} R[K])
$$
where $g$ is a unique natural operator
$$
g:J^{k-2}\Cla_{\tau} \f M\ucar{\f M} J^{k-1} \Lin \f E\ucar{\f M}
C^{(k-2,s-1)}_C\f M\ucar{\f M} C^{(k-1,r-1)}_L\f E
        \to   F\f E\,.\eqno{\END}
$$
\eTh

\smallskip

Further let us denote by $Z_i\f E$ the image of the operator
\begin{align*}
(\nab^{(i-2)}R[\Lam], & \nab^{(i-2)}R[K],\nab^{(i)}\Phi):
\\
J^{i-1}\Cla_{\tau}\f M\ucar{\f M}J^{i-1}\Lin\f E \ucar{\f M}
        J^i\f E^{p_1,p_2}_{q_1,q_2}
 & \to
C^{(i-2)}_C\f M\ucar{\f M} C^{(i-2)}_L\f E\ucar{\f M}
\f E^{p_1,p_2}_{q_1,q_2}\uten{\f M}\ten^iT^*\f M
\end{align*}
We shall put $Z^{(r)}\f E\byd Z_{0}\f E\ucar{\f M}\dots
\ucar{\f M} Z_{r}\f E$ and $\pr^r_k:Z^{(r)}\f E\to
Z^{(k)}\f E$, $r\ge k$, the canonical projection.

We denote by
$
J^{k-2}\Cla_{\tau} \f M\ucar{\f M} J^{k-2} \Lin \f E\ucar{\f M}
J^{k-1}\f E^{p_1,p_2}_{q_1,q_2}\ucar{\f M} Z^{(k,r)}\f E
$
the pullback of
$$
\pr^{r}_{k-1}:Z^{(r)}\f E\to Z^{(k-1)}\f E
$$
with respect to the surjective submersion, \cite{Jan04a},
$$
(\nab^{(k-3)}R[\Lam],\nab^{(k-3)}R[K],\nab^{(k-1)})
        :J^{k-2}\Cla_{\tau}\f M \ucar{\f M} J^{k-2}\Lin\f E
        \ucar{\f M}J^{k-1}\f E^{p_1,p_2}_{q_1,q_2}
        \to Z^{(k-1)}\f E\,.
$$
Then the second $k$-th order valued reduction theorem can be
formulated as follows, \cite{Jan04c},

\bTh\label{Th1.2}
Let $F$ be a $G$-gauge-natural bundle
of order $k\ge 1$ and let $r+1\ge k$. All
natural differential operators
$$
f:C^\infty(\Cla_{\tau}\f M\ucar{\f M} \Lin\f E\ucar{\f M}
        \f E^{p_1,p_2}_{q_1,q_2})\to C^\infty(F\f E)
$$
of order $r$ with respect sections of $\f E^{p_1,p_2}_{q_1,q_2}$
are of the form
\beq
f(j^{r-1}\Lam,j^{r-1}K,j^r\Phi) = g(j^{k-2}\Lam,
        j^{k-2}K, j^{k-1}\Phi,\nab^{(k-2,r-2)}
R[\Lam],\nab^{(k-2,r-2)}
R[K],
        \nab^{(k,r)} \Phi)\,
\eeq
where $g$ is a unique natural operator
$$
g:J^{k-2}\Cla_{\tau}\f M\ucar{\f M} J^{k-2}\Lin\f E \ucar{\f M}
J^{k-1}\f E^{p_1,p_2}_{q_1,q_2}\ucar{\f M}Z^{(k,r)}\f
M\to F\f E\,.\eqno{\END}
$$
\eTh

\bRm\label{Rm1.3}
The order $(r-1)$ of the above operators with respect to linear and
symmetric classical
connections is the minimal order we have to use. The second reduction
theorem can be easily generalized for any operators of orders $s_1$
or $s_2$
with respect to connections $\Lam$ or $K$,
respectively, where $s_1,s_2\ge r-1$, $s_1\ge s_2-1$. Then
$$
f(j^{s_1}\Lam,j^{s_2}K,j^r\Phi) =
        g(j^{k-2}\Lam, j^{k-2}K, j^{k-1}\Phi,\nab^{(k-2,s_1-1)}
R[\Lam], \nab^{(k-2,s_2-1)}
R[K],
        \nab^{(k,r)} \Phi)\,. \eqno{\END}
$$
\eRm

\bRm\label{Rm1.4}
The above higher order valued valued reduction
theorems deal with symmetric
classical connections on the base manifolds. If $\Lam$ is
a non-symmetric classical connection, then there is its unique splitting
$
  \Lam = \widetilde\Lam + T,
$
where $\widetilde\Lam$ is the symmetric classical connection
obtained from $\Lam$ by symmetrization, i.e.,
$\widetilde\Lam\col \mu\lam\nu = \tfr12\,(\Lam\col \mu\lam\nu +
\Lam\col \nu\lam\mu)$, and $T$ is the torsion $(1,2)$-tensor,
i.e., $T\col \mu\lam\nu = \tfr12\,(\Lam\col \mu\lam\nu -
\Lam\col \nu\lam\mu)$.
Then any finite order natural operator for
$\Lam$ and $K$ is of the form, $s\ge r-2$,
\begin{align*}
    f(j^s\Lam, j^rK) & = f(j^s\widetilde\Lam, j^rK, j^sT)=
\\
&=   g(j^{k-2}\widetilde\Lam, j^{k-2}K, j^{k-1} T,
        \widetilde\nab^{(k-2,s-1)} R[\widetilde\Lam],
        \widetilde\nab^{(k-2,r-1)} R[K],
        \widetilde\nab^{(k,s)} T)\,,\nonumber
\end{align*}
where $\widetilde{\ }$ refers to $\widetilde\Lam$.
{\ }\hfill\qedsymbol
\eRm

%--------------------------------------------------------------------%
\section{Natural classical connections on the total space of a
vector bundle}
\setcounter{equation}{0}
\setcounter{theorem}{0}
%--------------------------------------------------------------------%

A classical connection $D$ on $\f E$ is given by
$
  D:T\f E \to T^*\f E\ten TT\f E
$
over the identity of $T\f E$. In coordinates
\bEq
  D= d^C\ten (\der_C + D\col BAC\, \dot x^B\, \dot\der_A)\,,
        \quad D\col BAC\in C^\infty(\f E,\Rn)\,.
\eEq

Given a general linear connection $K$ on $\f E$ and
a classical connection $\Lam$ on the base manifol $\f M$
we have an induced natural
classical connection $D(\Lam,K)$ on $\f E$ given by,
\cite{Gan85, KolMicSlo93},

\bPr\label{Pr2.1}
There exists a unique classical connection $D=D(\Lam,K)$ on
the total space $\f E$ with the following properties
\begin{align*}
\nab^D_{h^K(X)} h^K(Y) = h^K(\nab^\Lam_X Y),& \quad
\nab^D_{h^K(X)} s^V = (\nab^K_X s)^V,
\\
\nab^D_{s^V} h^K(X) = 0,& \quad
\nab^D_{s^V} \sig^V = 0\,,
\end{align*}
for all vector fields $X,Y$ on $\f M$ and all sections $s,\sig$ of
$\f E$, where $h^K$ is the horisontal lift with respect to $K$,
$\nab^K, \nab^\Lam, \nab^D$ are covariant differentials with respect
to $K,\,\Lam\,, D$, respectively, and $s^V,\sig^V$ denote
the vertical lifts of the sections $s,\sig$, respectively.

In coordinates
\begin{align*}
D\col \mu\lam\nu & = \Lam \col\mu\lam\nu\,,
\quad
D\col \mu\lam{k}  = 0\,,
\quad
D\col j\lam\nu  = 0\,,
\quad
D\col j\lam{k}  = 0\,,
\\
 D\col\mu{i}\nu  & = \big(\der_\nu\, K\col pi\mu
        - K\col ri\nu\, K\col pr\mu
        + K\col pi\rho\,\Lam\col \mu\rho\nu\big)\, y^p\,,
\\
D\col \mu{i}k & = K\col ki\mu\,,
\quad
D\col j{i}\nu  = K\col ji\nu\,,
\quad
D\col j{i}k  = 0\,.
\end{align*}
\vglue-1.3\baselineskip{\ }\hfill\qedsymbol
\ePr

\bRm\label{Rm2.2}
$\Cla\f E$ is a $G$-gauge-natural bundle of order (2,2) and
$D(\Lam,K)$ defines the natural operator
$D:C^\infty(\Cla\f M\ucar{\f M}\Lin\f E)
        \to C^\infty(\Cla \f E)$
which is of order zero with respect to $\Lam$ and of order
one with respect to $K$.
\hfill\qedsymbol
\eRm

\medskip
The difference of any two classical connections
on $\f E$ is a tensor field on $\f E$ of the type $(1,2)$.
So, having the connection $D(\Lam,K)$,
all classical connections on $\f E$ naturally given by
$K$ and $\Lam$ are
of the type $D(\Lam,K)+ \Phi(\Lam,K)$,
where $\Phi(\Lam,K)$ is a natural (1,2)-tensor field on $\f E$.
Hence, the problem of classiffication of natural classical connections
on $\f E$
is reduced to the problem of classifffication of natural
tensor fields on $\f E$.

Any tensor field on $\f E$ is a section of a $G$-gauge-natural bundle
of order $(1,1)$. Then, by Theorem
\ref{Th1.2} and Remark \ref{Rm1.4}, we get

\bCr\label{Cr2.3}
Let $\Phi$ be a tensor field on $\f E$ naturally given by
a classical connection $\Lam$ on $\f M$
(in order $s$) and by
a  general linear connection $K$ on $\f E$ (in order $r$,
$s\ge r-2$).
Then
\beq
\Phi(u,j^s\Lam,j^rK) = \Psi (u,\widetilde\nab^{(s)} T,
\widetilde\nab^{(s-1)}R[\widetilde\Lam], K,
\widetilde\nab^{(r-1)} R[K])\,,
\eeq
where $u\in \f E$ and $\widetilde{\ }$ refers to the
classical symmetrized connection $\widetilde\Lam$.
\hfill\qedsymbol
\eCr

Now we can use the above Corollary \ref{Cr2.3} to classify
(1,2)-tensor fields on $\f E$. We have

\bLm\label{Lm2.4}
All (1,2)-tensor fields on $\f E$ naturally given by
$\Lam$ (in order $s$) and by
$K$ (in order $r$, $s\ge r-2$)
are of the maximal order 1 (with respect to both connections)
and form a 15-parameter
family of operators with the coordinate expression given by
\begin{align*}
\Phi(\Lam,K) & = \big(a_1\, T\col\mu\lam\nu +
        a_2\,\del^\lam_\mu T\col\rho\rho\nu +
        a_3\,\del^\lam_\nu T\col\mu\rho\rho\big)
        \,d^\mu\ten\der_\lam\ten d^\nu
\\
 & \quad +\bigg( y^i\,(
        b_{1}\, T\col\rho\rho\mu\, T\col \sig\sig\nu +
        b_{2}\, T\col\sig\rho\mu\, T\col \rho\sig\nu +
        b_{3}\, T\col\rho\rho\sig\, T\col \mu\sig\nu +
        c_{1}\, T\col\rho\rho{\mu;\nu} +
        c_{2}\, T\col\rho\rho{\nu;\mu}
\\
&\quad
        + c_{3}\, T\col\mu\rho{\nu;\rho} +
        d_1\, \widetilde R\col\rho\rho{\mu\nu}
        + d_2\, \widetilde R\col\mu\rho{\rho\nu} +
        e_2\, R\col pp{\mu\nu}) +
        e_1\, R\col ji{\mu\nu}\, y^j
\\
&\quad  +
        (a_3- h_{2})\, T\col \mu\rho\rho\, K\col ji\nu\, y^j +
        (a_2- h_{1})\, T\col \rho\rho\nu\, K\col ji\mu\, y^j
        + a_1\, T\col \mu\rho\nu\, K\col ji\rho\, y^j\bigg)
        \,d^\mu\ten\der_i\ten d^\nu
\\
        &\quad+ h_{1}\, \del^i_j\, T\col\rho\rho\nu
        \,d^j\ten\der_i\ten d^\nu
        + h_{2}\, \del^i_k\, T\col\rho\rho\mu
        \,d^\mu\ten\der_i\ten d^k
\,,
\end{align*}
where $a_i, b_i, c_i$, $i=1,2,3$, $d_j,e_j,h_j$, $j=1,2$,
are real coefficients.
\eLm

\bPf
By the general theory of natural operators, \cite{KolMicSlo93},
we have to classify all equivariant mappings between standard
fibers of $G$-gauge-natural bundles in question.
$T^*\f E\ten T\f E\ten T^*\f E$ is a $G$-gauge-natural bundle
of order $(1,1)$ where the action of the group $W^{(1,1)}_{m,n} G$
on the standard fiber $S_F=\Rn^{(m+n)*}\ten\Rn^{(m+n)}\ten\Rn^{(m+n)*}$
is given by
\begin{align}
\ba \Phi_j{}^i{}_k & = a^i_r\, \Phi_s{}^r{}_t\,
        \tilde a^s_j\, \tilde a^t_k
        + a^i_{p\rho}\, y^p\, \Phi_s{}^\rho{}_t\,
        \tilde a^s_j\, \tilde a^t_k\,,\label{Eq2.1}
\\
\ba \Phi_j{}^i{}_\nu & = a^i_r\, \Phi_s{}^r{}_\tau\,
        \tilde a^s_j\, \tilde a^\tau_\nu
        + a^i_r\, \Phi_s{}^r{}_t\, \label{Eq2.2}
        \tilde a^s_j\, \tilde a^t_{p\nu}\, a^p_k\, y^k
\\
        & \quad + a^i_{p\rho}\, y^p\, \Phi_s{}^\rho{}_t\,
        \tilde a^s_j\, \tilde a^t_{k\nu}\, a^k_q\, y^q
        + a^i_{p\rho}\, y^p\, \Phi_s{}^\rho{}_\tau\,
        \tilde a^s_j\, \tilde a^\tau_{\nu}\,, \nonumber
\\
\ba \Phi_\mu{}^i{}_k & = a^i_r\, \Phi_\sig{}^r{}_t\,
        \tilde a^\sig_\mu\, \tilde a^t_k
        + a^i_r\, \Phi_s{}^r{}_t\,           \label{Eq2.3}
        \tilde a^s_{p\mu}\, a^p_q\, y^q \tilde a^t_k
\\
      &\quad  + a^i_{p\rho}\, y^p\, \Phi_s{}^\rho{}_t\,
        \tilde a^s_{q\mu}\,a^q_l\, y^l\, \tilde a^t_{k}
        + a^i_{p\rho}\, y^p\, \Phi_\sig{}^\rho{}_t\,
        \tilde a^\sig_\mu\, \tilde a^t_{k}\,,\nonumber
\\
\ba \Phi_\mu{}^i{}_\nu & = a^i_r\, \Phi_\sig{}^r{}_\tau\,
        \tilde a^\sig_\mu\, \tilde a^\tau_\nu
        + a^i_r\, \Phi_s{}^r{}_\tau\,\label{Eq2.4}
        \tilde a^s_{p\mu}\, a^p_q\, y^q \tilde a^\tau_\nu
\\
   &\quad        + a^i_r\, \Phi_\sig{}^r{}_t\,
        \tilde a^\sig_{\mu}\, \tilde a^t_{p\nu}\, a^p_l\, y^l
        + a^i_r\, \Phi_s{}^r{}_t\,   \nonumber
        \tilde a^s_{p\mu}\,a^p_l\, y^l \tilde a^t_{q\nu}\, a^p_m\, y^m
\\
    &\quad        + a^i_{p\rho}\, y^p\, \Phi_s{}^\rho{}_t\,
        \tilde a^s_{q\mu}\,a^q_l\, y^l\, \tilde a^t_{n\nu}\, a^n_m\, y^m
        + a^i_{p\rho}\, y^p\, \Phi_\sig{}^\rho{}_t\,\nonumber
        \tilde a^\sig_\mu\, \tilde a^t_{q\nu}\, a^q_l\, y^l
\\
     &\quad   + a^i_{p\rho}\, y^p\, \Phi_s{}^\rho{}_\tau\,
        \tilde a^s_{q\mu}\,a^q_l \, y^l \tilde a^\tau_{\nu}
        + a^i_{p\rho}\, y^p\, \Phi_\sig{}^\rho{}_\tau\,
        \tilde a^\sig_{\mu}\, \tilde a^\tau_{\nu}\,,\nonumber
\\
\ba \Phi_j{}^\lam{}_k & = a^\lam_\rho\, \Phi_s{}^\rho{}_t\,
        \tilde a^s_j\, \tilde a^t_k\,,\label{Eq2.5}
\\
\ba \Phi_j{}^\lam{}_\nu & = a^\lam_\rho\, \Phi_s{}^\rho{}_\tau\,
        \tilde a^s_j\, \tilde a^\tau_\nu
        + a^\lam_\rho\, \Phi_s{}^\rho{}_t\,
        \tilde a^s_j\, \tilde a^t_{p\nu}\, a^p_k\, y^k\,,\label{Eq2.6}
\\
\ba \Phi_\mu{}^\lam{}_k & = a^\lam_\rho\, \Phi_\sig{}^\rho{}_t\,
        \tilde a^\sig_\mu\, \tilde a^t_k
        + a^\lam_\rho\, \Phi_s{}^\rho{}_t\,
        \tilde a^s_{p\mu}\, a^p_q\, y^q \tilde a^t_k\,, \label{Eq2.7}
\\
\ba \Phi_\mu{}^\lam{}_\nu & = a^\lam_\rho\, \Phi_\sig{}^\rho{}_\tau\,
        \tilde a^\sig_\mu\, \tilde a^\tau_\nu
        + a^\lam_\rho\, \Phi_s{}^\rho{}_\tau\, \label{Eq2.8}
        \tilde a^s_{p\mu}\, a^p_q\, y^q \tilde a^\tau_\nu
\\
    & \quad   + a^\lam_\rho\, \Phi_\sig{}^\rho{}_t\,
        \tilde a^\sig_{\mu}\, \tilde a^t_{p\nu}\, a^p_l\, y^l
        + a^\lam_\rho\, \Phi_s{}^\rho{}_t\,
        \tilde a^s_{p\mu}\,a^p_l\, y^l \tilde a^t_{q\nu}\,
        a^q_m\, y^m\,.
          \nonumber
\end{align}

By Corollary \ref{Cr2.3}  equivariant
mappings $\Phi:\Rn^n\car T^s_m Q\car T^r_m R\to S_F$
between standard fibers
are in the form
$$
\Phi\col BAC (y^i,\Lam \col \mu\lam{\nu,\f{\alp}},
        K\col ji{\mu,\f{\bet}})
        =\Psi \col BAC(y^i,
        T\col \mu\lam{\nu;\f{\gam}},
        \widetilde R\col \mu\lam{\nu\kappa;\f{\del}},
        K\col ji\mu,R\col ji{\lam\mu;\f{\eps}})\,,
$$
where $\f{\alp}, \f{\bet},\f{\gam}, \f{\del}, \f{\eps}$
are multiindices of the ranks
$0\le\, \parallel\!\!\f{\alp}\!\!\parallel,
\parallel\!\!\f{\gam}\!\!\parallel\, \le s$,
$0\le\, \parallel\!\!\f{\bet}\!\!\parallel\,\le  r$,
$0\le\, \parallel\!\!\f{\del}\!\!\parallel\,\le  s-1$,
$0\le\, \parallel\!\!\f{\eps}\!\!\parallel\,\le  r-1$, and
$(\Lam \col \mu\lam{\nu,\f{\alp}})$ or
$(K\col ji{\mu,\f{\bet}})$ are the induced coordinates on
$T^s_mQ$ or $T^r_m R$, respectively.

First, let us consider the equivariance of $\Psi\col BAC$ with
respect to the fiber homotheties, i.e.,
with respect to $(c\, \del^i_j)$.
Then, from
(\ref{Eq2.5}), we get
$$
 c^{-2}\, \Psi\col j\lam{k} =\Psi \col j\lam{k}(c\, y^i,
        T\col \mu\lam{\nu;\f{\gam}},
        \widetilde R\col \mu\lam{\nu\kappa;\f{\del}},
        K\col ji\mu,R\col ji{\lam\mu;\f{\eps}})\,.
$$
Multiplying it by $c^2$ and letting $c\to 0$, we obtain
\bEq\label{Eq2.9}
\Psi\col j\lam{k} = \Phi\col j\lam{k} = 0\,.
\eEq

Similarly, from (\ref{Eq2.6}), (\ref{Eq2.7}) and (\ref{Eq2.1}), we have
\begin{align*}
 c^{-1}\, \Psi\col j\lam{\nu} & =\Psi \col j\lam{\nu}(c\, y^i,
        T\col \mu\lam{\nu;\f{\gam}},
        \widetilde R\col \mu\lam{\nu\kappa;\f{\del}},
        K\col ji\mu,R\col ji{\lam\mu;\f{\eps}})\,,
\\
 c^{-1}\, \Psi\col \mu\lam{k} & =\Psi \col \mu\lam{k}(c\, y^i,
        T\col \mu\lam{\nu;\f{\gam}},
        \widetilde R\col \mu\lam{\nu\kappa;\f{\del}},
        K\col ji\mu,R\col ji{\lam\mu;\f{\eps}})\,,
\\
 c^{-1}\, \Psi\col jik & =\Psi \col jik (c\, y^i,
        T\col \mu\lam{\nu;\f{\gam}},
        \widetilde R\col \mu\lam{\nu\kappa;\f{\del}},
        K\col ji\mu,R\col ji{\lam\mu;\f{\eps}})\,,
\end{align*}
which, by multiplying by $c$ and letting $c\to 0$, imply
\begin{align}\label{Eq2.10}
\Psi\col j\lam{\nu} & = \Phi\col j\lam{\nu} = 0\,,
\quad
\Psi\col \mu\lam{k}  = \Phi\col \mu\lam{k} = 0\,,
\quad
\Psi\col jik = \Phi\col jik = 0\,.
\end{align}

Further, from (\ref{Eq2.2}), (\ref{Eq2.3}) and (\ref{Eq2.8}), we have
\begin{align*}
\Psi\col ji{\nu} & =\Psi \col ji{\nu}(c\, y^i,
        T\col \mu\lam{\nu;\f{\gam}},
        \widetilde R\col \mu\lam{\nu\kappa;\f{\del}},
        K\col ji\mu,R\col ji{\lam\mu;\f{\eps}})\,,
\\
\Psi\col \mu{i}{k} & =\Psi \col \mu{i}{k}(c\, y^i,
        T\col \mu\lam{\nu;\f{\gam}},
        \widetilde R\col \mu\lam{\nu\kappa;\f{\del}},
        K\col ji\mu,R\col ji{\lam\mu;\f{\eps}})\,,
\\
\Psi\col \mu\lam\nu & =\Psi \col \mu\lam\nu (c\, y^i,
        T\col \mu\lam{\nu;\f{\gam}},
        \widetilde R\col \mu\lam{\nu\kappa;\f{\del}},
        K\col ji\mu,R\col ji{\lam\mu;\f{\eps}})\,,
\end{align*}
which implies, by letting $c\to 0$, that $\Psi\col ji{\nu},
\Psi\col \mu{i}{k}$ and $\Psi\col \mu\lam\nu$ are independent of $y^i$.

Finally, from (\ref{Eq2.4}), we have
$$
c\,\Psi\col \mu{i}\nu  =\Psi \col \mu{i}\nu (c\, y^i,
        T\col \mu\lam{\nu;\f{\gam}},
        \widetilde R\col \mu\lam{\nu\kappa;\f{\del}},
        K\col ji\mu,R\col ji{\lam\mu;\f{\eps}})\,.
$$
By the homogeneous function theorem, \cite{KolMicSlo93},
$\Psi\col \mu{i}\nu$ is linear in $y^p$, i.e.,
$$
\Psi\col \mu{i}\nu  =F \col \mu{i}{\nu p} (
        T\col \mu\lam{\nu;\f{\gam}},
        \widetilde R\col \mu\lam{\nu\kappa;\f{\del}},
        K\col ji\mu,R\col ji{\lam\mu;\f{\eps}})\, y^p\,.
$$
So we have equivariant functions
$\Psi\col ji{\nu}$,
$\Psi\col \mu{i}{k}$, $\Psi\col \mu\lam\nu$ and $F \col \mu{i}{\nu p}$
of variables
$T\col \mu\lam{\nu;\f{\gam}}$,
\linebreak
$\widetilde R\col \mu\lam{\nu\kappa;\f{\del}}$,
        $K\col ji\mu$, $R\col ji{\lam\mu;\f{\eps}}$.

Now we consider the base homotheties, i.e., the equivariance with respect
to $(c\, \del^\lam_\mu)$.
Then
$$
c^{-1}\,\Psi\col \mu\lam\nu  =\Psi \col \mu\lam\nu (
        c^{-(\parallel\f{\gam}\parallel+1)}\,T\col \mu\lam{\nu;\f{\gam}},
        c^{-(\parallel\f{\del}\parallel+2)}\,
        \widetilde R\col \mu\lam{\nu\kappa;\f{\del}},
        c^{-1}\, K\col ji\mu,
        c^{-(\parallel\f{\eps}\parallel+2)}\,R\col ji{\lam\mu;\f{\eps}}
        )\,.
$$
The homogeneous function theorem implies that $\Psi\col \mu\lam\nu$
is a polynomial function with exponents
$a_i$ in
$T\col \mu\lam{\nu;\f{\gam}}$,
$i= {\parallel\!\!\f{\gam}\!\!\parallel}$,
$b_j$ in
$R\col \mu\lam{\nu\kappa;\f{\del}}$,
$j={\parallel\!\!\f{\del}\!\!\parallel}$,
$c$ in $K\col ji\mu$
and $d_k$ in
$R\col ji{\lam\mu;\f{\eps}}$,
$k={\parallel\!\!\f{\eps}\!\!\parallel}$,
such that
\bEq\label{Eq2.11}
-1 = - \sum^{s}_{i=0}
        (i+1)\,a_i
        - \sum^{s-1}_{j=0}
        (j+2)\, b_j
        - c -
         \sum^{r-1}_{k=0}
        (k+2)\,
        d_k \,.
\eEq
The above equation (\ref{Eq2.11}) has solutions in integers only for
$a_0=1$ and the other coefficients vahish or $c=1$ and the other
coefficients vanish. This implies that
$\Psi\col \mu\lam\nu$ is of the type
$$
\Psi\col \mu\lam\nu = A^{\lam\sig\tau}_{\mu\nu\rho}\,
T\col \sig\rho\tau + L^{\lam j\rho}_{\mu\nu i}\, K\col ji\rho\,,
$$
where $A^{\lam\sig\tau}_{\mu\nu\rho}$ and $L^{\lam j\rho}_{\mu\nu i}$
are absolute invariant tensors, i.e., products of the Kronecker symbols.
Then
$$
\Psi\col \mu\lam\nu = a_1\, T\col \mu\lam\nu
        + a_2\, \del^\lam_\mu\, T\col \rho\rho\nu
        + a_3\, \del^\lam_\nu\, T\col \rho\rho\mu
        + l_1\,\del^\lam_\mu K\col ii\nu
        + l_2\,\del^\lam_\nu K\col ii\mu\,.
$$
Finally, the equivariance of $\Psi\col \mu\lam\nu$ with respect to
elements of the type $(\del^i_j, \del ^\lam_\mu, a^i_{j\lam})$
implies
$l_1=l_2=0$ and we have
\bEq\label{Eq2.12}
\Psi\col \mu\lam\nu = \Phi\col \mu\lam\nu =  a_1\, T\col \mu\lam\nu
        + a_2\, \del^\lam_\mu\, T\col \rho\rho\nu
        + a_3\, \del^\lam_\nu\, T\col \rho\rho\mu\,.
\eEq

Similarly, for $\Psi\col ji\nu$ and $\Psi\col \mu{i}{k}$,
we get polynomial functions with exponents satisfying the same equation
(\ref{Eq2.11}). So we have
\begin{align*}
\Psi\col ji\nu & = H^{i\sig\tau}_{j\nu\rho}\,
T\col \sig\rho\tau + M^{im\rho}_{j\nu k}\, K\col mk\rho\,,
\quad
\Psi\col \mu{i}k  = H^{i\sig\tau}_{\mu k\rho}\,
T\col \sig\rho\tau + M^{im\rho}_{\mu kp}\, K\col mp\rho\,,
\end{align*}
where coefficients are absolute invariant tensors, i.e.,
\begin{align*}
\Psi\col ji\nu & = h_{1}\,\del^i_j\,
T\col \rho\rho\nu + m_1\, K\col ji\nu + m_2\,\del^i_j\, K\col pp\nu
\,,
\\
\Psi\col \mu{i}k & = h_{2}\,\del^i_k\,
T\col \rho\rho\mu + m_3\, K\col ki\mu + m_4\,\del^i_k\, K\col pp\mu
\,.
\end{align*}
The equivariance with respect to
elements of the type $(\del^i_j, \del ^\lam_\mu, a^i_{j\lam})$
implies
$m_i=0$, $i=1,\dots,4$, and we have
\bEq\label{Eq2.13}
\Psi\col ji\nu = \Phi\col ji\nu  = h_{1}\,\del^i_j\,
T\col \rho\rho\nu\,,\quad
\Psi\col \mu{i}k = \Phi\col \mu{i}k = h_{2}\,\del^i_k\,
T\col \rho\rho\mu\,.
\eEq

Finally, let us discuss $F\col \mu{i}{\nu j}$. Considering
base homotheties and the homogeneous function theorem we get that
$F\col \mu{i}{\nu j}$ is polynomial with exponents satisfying
(\ref{Eq2.11}) with $-2$ on the left hand side.
%$$
%-2 = - \sum^{s}_{i=0} (i+1)\, a_i - \sum^{s-1}_{j=0} (j+2)\,b_j
%        - c -
%         \sum^{r-1}_{k=0} (k+2)\, d_k \,.
%$$
Then we have the following 6 possible solutions:
$a_0=2$ and the other exponents vanish;
$a_1=1$ and the other exponents vanish;
$a_0=1$, $c=1$ and the other exponents vanish;
$b_0=1$ and the other exponents vanish;
$c=2$ and the other exponents vanish;
$d_0=1$ and the other exponents vanish.
Then
\begin{align*}
F\col \mu{i}{\nu j} & =
        B^{i\sig_1\tau_1\sig_2\tau_2}_{\mu\nu j\rho_1\rho_2}
        \, T\col{\sig_1}{\rho_1}{\tau_1}\,T\col{\sig_2}{\rho_2}{\tau_2}
        +C^{i\sig\tau_1\tau_2}_{\mu\nu j\rho}\,
        T\col\sig\rho{\tau_1;\tau_2}
        + N^{i\sig \tau_1 q\tau_2}_{\mu\nu j\rho p}\,
        T\col \sig\rho{\tau_1}\, K\col qp{\tau_2}
\\
    &\quad +
        D^{i\sig\tau_1\tau_2}_{\mu\nu j \rho}\,
        \tilde R\col\sig\rho{\tau_1\tau_2}
        +  P^{iq_1\tau_1 q_2\tau_2}_{\mu\nu j p_1 p_1}\,
        K\col {q_1}{p_1}{\tau_1}\,K\col {q_2}{p_2}{\tau_2} +
        E^{iq\tau_1\tau_2}_{\mu\nu jp}\,
        R\col qp{\tau_1\tau_2}\,,
\end{align*}
where all coefficients are absolute invariant tensors, i.e.,
\begin{align*}
F\col \mu{i}{\nu j} & = \del^i_j\,(
        b_{1}\, T\col \rho\rho\mu\, T\col \sig\sig\nu
        +b_{2}\, T\col \sig\rho\mu\, T\col \rho\sig\nu
        +b_{3}\, T\col \rho\rho\sig\, T\col \mu\sig\nu)
\\
     &\quad +
        \del^i_j\,(c_{1}\, T\col\rho\rho{\mu;\nu}
        +c_{2}\, T\col\rho\rho{\nu;\mu} +
        c_{3}\, T\col\mu\rho{\nu;\rho})
\\
     &\quad +
        \del^i_j\,(n_{1}\, T\col\rho\rho{\mu}\, K\col pp\nu
        + n_{2}\, T\col\rho\rho{\nu}\, K\col pp\mu
        + n_{3}\, T\col\mu\rho{\nu}\, K\col pp\rho)
\\
   & \quad +
        n_{4}\, T\col\mu\rho{\rho}\, K\col ji\nu
        + n_{5}\, T\col\rho\rho{\nu}\, K\col ji\mu
        + n_{6}\, T\col\mu\rho{\nu}\, K\col ji\rho
\\
  & \quad +
        d_1\, \del^i_j\,\widetilde R\col \rho\rho{\mu\nu} +
        d_2\, \del^i_j\,\widetilde R\col \mu\rho{\rho\nu}+
        e_1\, \del^i_j\, R\col pp{\mu\nu}+
        e_2\, R\col ji{\mu\nu}
\\
   & \quad +
        \del^i_j\, ( p_1\, K\col pp\mu\, K\col qq\nu +
        p_2\, K\col qp\mu\, K\col pq\nu)
\\
   & \quad +
        p_3\, K\col ji\mu\, K\col pp\nu +
        p_{4}\, K\col ji\nu\, K\col pp\mu +
        p_{5}\, K\col pi\mu\, K\col jp\nu +
        p_{6}\, K\col pi\nu\, K\col jp\mu\,.
\end{align*}

The equivariance with respect to
elements of the type $(\del^i_j, \del ^\mu_\nu, a^i_{j\mu})$
implies $n_1=n_2=n_3=0,\,
n_4 = a_3-h_{2}, \, n_5 = a_2- h_{1},\, n_6= a_1, p_i =0$,
$i=1\dots,6$, and the other coefficients are arbitrary. Then
\begin{align}\label{Eq2.14}
\Psi \col \mu{i}\nu & =\Phi \col \mu{i}\nu=
        F\col \mu{i}{\nu j}\,y^j  =
        y^i\,(
        b_{1}\, T\col \rho\rho\mu\, T\col \sig\sig\nu
        +b_{2}\, T\col \sig\rho\mu\, T\col \rho\sig\nu
        +b_{3}\, T\col \rho\rho\sig\, T\col \mu\sig\nu)
\\
     &\quad +
        y^i\,(c_{1}\, T\col\rho\rho{\mu;\nu}
        +c_{2}\, T\col\rho\rho{\nu;\mu}
        + c_{3}\, T\col\mu\rho{\nu;\rho})
        \nonumber
\\
   & \quad +
        y^j\,\big((a_{3}- h_{2})\, T\col\mu\rho{\rho}\, K\col ji\nu
        + (a_{2}- h_{1})\, T\col\rho\rho{\nu}\, K\col ji\mu +
        a_{1}\, T\col\mu\rho{\nu}\, K\col ji\rho\big)
        \nonumber
\\
  & \quad +
        d_1\, y^i\,\widetilde R\col \rho\rho{\mu\nu} +
        d_2\, y^i\,\widetilde R\col \mu\rho{\rho\nu}+
        e_1\, y^i\, R\col pp{\mu\nu} +
        e_2\, R\col ji{\mu\nu}\, y^j \,.
        \nonumber
\end{align}

Summerizing all results (\ref{Eq2.9}), (\ref{Eq2.10}),
(\ref{Eq2.12}), (\ref{Eq2.13}) and  (\ref{Eq2.14})
we get Lemma \ref{Lm2.4}.
\ePf

As a direct consequence of Lemma \ref{Lm2.4} we have

\bCr\label{Cr2.5}
All natural operators transforming $\Lam$ and $K$
into classical connections on $\f E$
are of the maximal order one
and form 15-parameter family
$$
  \widetilde{D}(\Lam, K)=D(\Lam,K) + \Phi(\Lam,K)\,,
$$
where $D(\Lam,K)$ is the connection given by Proposition \ref{Pr2.1}
and $\Phi(\Lam,K)$ is the 15-parameter family of natural tensor
fields given by Lemma \ref{Lm2.4}. \hfill \END
\eCr

\bRm\label{Rm2.6}
In \cite{GanKol91} (see also \cite{KolMicSlo93}, Proposition 54.3)
the same result was obtained
by direct calculations without using the reduction theorems.
Our result coincides with the result of \cite{GanKol91, KolMicSlo93}
but our base of the 15-parameter family of
operators of Lemma \ref{Lm2.4} differ from the base used
in \cite{GanKol91, KolMicSlo93}.
\hfill\qedsymbol
\eRm

Now we shall describe the geometrical construction of the operators
from Lemma \ref{Lm2.4}. First let us recall that we have the
canonical immersions $\iota_{T^*\f M}:T^*\f M\to T^*\f E$ and
$\iota_{V\f E}:V\f E\to T\f E$. The connection $K$ defines
the horisontal lift
$
h^K:T\f M\to T\f E
$
and the vertical projection
$
\nu_K:T\f E\to V\f E\,.
$

Then we have the following subfamilies of natural operators given by
$\Lam$ and $K$.

A) $\Lam$ gives 3-parameter family of (1,2)-tensor fields on $\f M$,
\cite{KolMicSlo93}, given by
\bEq
S(\Lam) = a_1\, T + a_2 I_{T\f M}\ten \hat T
        + a_3\, \hat T\ten I_{T\f M}\,,
\eEq
where $T$ is the torsion tensor of $\Lam$, $\hat T$ is its contraction
and $I_{T\f M}: M\to T\f M\ten T^*\f M$ is the identity tensor.
Then the first 3-parameter subfamily of operators
of Lemma \ref{Lm2.4}
is given by
$
h^K
(S(\Lam))\equiv\iota_{T^*\f M}\ten h^K\ten \iota_{T^*\f M}
(S(\Lam))$.

B) $\Lam$ and $K$ define naturally
the following 9-parameter
family of (0,2) tensor fields on $\f M$, \cite{KolMicSlo93},
given by
\begin{align*}
G(\Lam, K) & = b_{1}\, C^{12}_{13} (T\ten T) +
        b_{2}\, C^{12}_{31} (T\ten T) +
        b_{3}\, C^{12}_{12} (T\ten T)
\\
 &\quad +
        c_{1}\, C^{1}_{1} \widetilde\nab T +
        c_{2}\, \overline{C^{1}_{1}  \widetilde\nab T} +
        c_{3}\, C^{1}_{3} \widetilde\nab  T +
        d_1\, C^1_1 R[\widetilde\Lam] + d_2\, C^1_2 R[\widetilde\Lam]+
        e_1\, C^1_1 R[K]\,,
\end{align*}
where $C^{ij}_{kl}$ is the contraction with respect to indicated indices
and $\overline{C^{1}_{1}  \widetilde\nab T}$
denotes the conjugated tensor obtained
by the exchange of subindices.
The second 9-parameter subfamily of operators
from Lemma \ref{Lm2.4} is then
given by
$
L\ten G(\Lam,K)\equiv
\iota_{T^*\f M}\ten \iota_{V\f E}\ten \iota_{T^*\f M}
(L\ten G(\Lam,K))$,
where $L=y^i\, \der_i$ is the Liouville
vertical vector field on $\f E$.

C) The value of the curvature tensor $R[K]$ applied on the Liouville
vector field is in $T^*\f M\ten V\f E\ten T^*\f M$.
Then $R[K](L)\equiv\iota_{T^*\f M}\ten \iota_{V\f E}\ten \iota_{T^*\f M}
(R[K](L))$ is the operator standing by $e_2$ in
Lemma \ref{Lm2.4}.

D) Finally, if we consider $\nu_K$ as the vertical valued 1-form
$
  \nu_K:\f E\to T^*\f E\ten V\f E
$
with coordinate expression
$$
  \nu_K= (d^i - K\col ji\lam \,y^j\, d^\lam)\ten\der_i\,,
$$
the last 2-parameter subfamily of operators from Lemma \ref{Lm2.4}
is obtained by applying the morphism
$
\iota_{T^*\f M}\ten\iota_{V\f E}\ten \id_{T^*\f E}
$
on
$$
H(\Lam, K)= h_{1} \,\nu_K\ten \hat T + h_{2}\, \hat T\ten \nu_K
      \,.
$$

Summarizing the above constructions we get

\bTh\label{Th2.7}
All classical connections on $\f E$ naturally given by
$\Lam$ (in order $s$) and by
$K$ (in order $r$,
$s\ge r-2$) are of the maximal order one and are of
the form
\begin{align*}
\widetilde{D}(\Lam,K)=D(\Lam,K) +
        h^K\big(S(\Lam)\big)
 +L\ten G(\Lam,K) + e_2\, R[K](L)
       + H(\Lam,K)\,.
\end{align*}
\vglue-1.3\baselineskip{\ }\hfill\END
\eTh

\bCr\label{Cr2.8}
All natural operators transforming a  general linear
connection $K$ on $\f E$
and a symmetric classical connection $\Lam$ on
$\f M$ into classical connections on $\f E$
are of the maximal order one
and form the following 4-parameter family
$$
\widetilde{D}(\Lam,K)
        =  D(\Lam,K) + L\ten\big(d_1\, C^1_1 R[\widetilde\Lam]
        + d_2\,  C^1_2 R[\widetilde\Lam]+
         e_1\, C^1_1 R[K]\big) + e_2\, R[K](L)\,.\eqno{\qedsymbol}
$$
\eCr

%--------------------------------------------------------------------%
\section{Natural connections on the 1st jet prolongation of vector
bundles}
\setcounter{equation}{0}
\setcounter{theorem}{0}
%--------------------------------------------------------------------%

Assume the 1-jet prolongation $\pi^1_0:J^1\f E \to \f E$ with
the induced fiber coordinate chart $(x^\lam, y^i; y^i_\lam)$.
Then $J^1\f E$ is a (1,1)-order $G$-gauge-natural bundle with the
standard fiber $\Rn^n\car \Rn^n\ten\Rn^{m*}$ and
the induced action of $W^{(1,1)}_m G$ given in coordinates by
$$
\bar y^i= a^i_p\, y^p\,,\quad
        \bar y^i_\lam = a^i_p\, y^p_\rho\, \tilde a^\rho_\lam +
        a^i_{p\rho}\, y^p\, \tilde a^\rho_\lam\,.
$$

A connection $\Gam$ on $J^1 \f E$ is given by
$
  \Gam: J^1\f E\to T^*\f E\ten TJ^1\f E
$
over the identity of $T\f E$ with coordinate expression
\bEq
  \Gam = d^A\ten \big( \der_A + \Gam_A{}_\lam^i\,\der^\lam_i\big)
\,,\quad \Gam_A{}^i_\lam\in C^\infty(J^1\f E,\Rn)\,.
\eEq
A connection $\Gam$ is affine if and only if
$
  \Gam_A{}^i_\lam = \Gam_A{}^i_\lam{}^\mu_p\, y^p_\mu +
             \Gam_A{}^i_\lam{}^._.\,.
$

\bRm\label{Rm3.1}
If $\Gam_1$ and $\Gam_2$ are two connections on $J^1\f E$, then
the difference
$
\Gam_1-\Gam_2:J^1\f E\to T^*\f E \ten VJ^1\f E
$
and, by using the identification
$
  VJ^1\f E = J^1\f E\ucar{\f E} T^*\f M\uten{\f E} V\f E\,,
$
we get
$$
\Gam_1-\Gam_2=\phi:J^1\f E\to T^*\f E\ten T^*\f M\ten V\f E\,.
\eqno{\qedsymbol}
$$
\eRm

In \cite{JanMod96} we have described a natural operator $\chi$
transforming a classical connection on the total space of
a fibered
manifold and a classical connection on the base manifold
into a connection on the 1st jet prolongation of the
fibered manifold. Applying this operator on a classical
connection $D$ on the total space of a vector bundle $\f E\to \f M$
we get

\bPr\label{Pr3.2}
Let $D$ be a classical connection on $\f E$ and $\Lam$ be
a classical connection on $\f M$. Then we have a connection
$\Gam(\Lam,D)=\chi(D)$ on $\pi^1_0:J^1\f E\to \f E$ such that
$$
  \Gam_A{}^i_\lam = D\col Aij\, y^j_\lam + D\col Ai\lam -
        y^i_\mu\, (D\col A\mu{j}\, y^j_\lam + D\col A\mu\lam)\,.
$$
\ePr

\bRm\label{Rm3.3}
The connection $\Gam(\Lam,D)$ is independent of $\Lam$, but the
geometric construction of the operator $\chi$ depends on $\Lam$
essentially, see \cite{JanMod96}.
\hfill\qedsymbol
\eRm

Now, applying the operator $\chi$ on the connection $D(\Lam,K)$
from Proposition \ref{Pr2.1}, we get

\bTh\label{Th3.4}
A  general linear connection $K$ on $\f E$ and
a classical connection
$\Lam$ on $\f M$ give naturally the connection
$\Gam(\Lam,K)=\chi(D(\Lam,K))$ on $J^1\f E$ with the coordinate
expression
\begin{align*}
\Gam_\mu{}^i_\lam & = K\col ji\mu\, y^j_\lam +
        (\der_\lam K\col ji\mu - K\col pi\lam\, K\col jp\mu
        + K\col ji\rho\, \Lam\col \mu\rho\lam)\, y^j\,,
\\
\Gam_j{}^i_\lam & = K\col ji\lam\,.
\end{align*}
\vglue-1.3\baselineskip
{\ }\hfill\qedsymbol
\eTh

\bCr\label{Cr3.5}
All natural operators transforming a (general)
linear connection $K$
on $\f E$ (in order $r$) and a classical connection
$\Lam$ on $\f M$ (in order $s$, $s\ge r-2$) into connections
on $J^1\f E$ are of the form
$$
\widetilde{\Gam}(\Lam,K) = \Gam(\Lam,K) + \phi(\Lam,K)\,,
$$
where $\Gam(K,\Lam)$ is the connection from Theorem \ref{Th3.4}
and $\phi(\Lam,K)$ is a
natural operator
$$
\phi: J^1\f E\ucar{\f M} J^s\Cla \f M\ucar{\f M} J^r\Lin \f E  \to
        T^*\f E\ten T^*\f M\ten V\f E\,.
        \eqno{\qedsymbol}
$$
\eCr

So to classify natural connections on $J^1\f E$ it is sufficient to
classify natural operators $\phi$.

\bLm\label{Lm3.6}
All natural tensor fields
$\phi(\Lam,K):J^1\f E\to T^*\f E\ten T^*\f M\ten V\f E$
naturally given by
a classical connection $\Lam$ on $\f M$ (in order $s$) and by
a  general linear connection $K$ on $\f E$ (in order $r$,
$s\ge r-2$)  are of the maximal
order one (with respect to both connections) and form a 14-parameter
family of operators with coordinate expression given by
\begin{align*}
\phi(\Lam,K) & = \bigg(\big(a_1\, T\col \lam\rho\mu
        +  a_2\, T\col \sig\sig\mu \, \del^\rho_\lam
        + a_3\, T\col \lam\sig\sig \, \del^\rho_\mu\big)\,y^i_\rho
        + y^i\, \big(b_1\,  T\col \rho\rho\lam \, T\col \sig\sig\mu
        + b_2\, T\col \sig\rho\lam \, T\col \rho\sig\mu
\\
&\quad         + b_3\,  T\col \rho\rho\sig \, T\col \lam\sig\mu
      + c_1\, T\col \rho\rho{\lam;\mu}
        + c_2\, T\col \rho\rho{\mu;\lam}
        + c_3\, T\col \lam\rho{\mu;\rho}\big)
 \\
&\quad        - \big(a_3\, T\col \rho\rho{\lam}\, K\col ji\mu
        + (a_2+h_{1})\, T\col \rho\rho{\mu}\, K\col ji\lam
        - a_1\, T\col \lam\rho{\mu}\, K\col ji\rho \big)\, y^j
\\
   &\quad + y^i\,\big(d_1 \, R\col \rho\rho{\lam\mu}
        + d_2 \, R\col \lam\rho{\rho\mu}\big)
        + e_1 \,y^i\, R\col pp{\lam\mu}
        + e_2 \, R\col ji{\lam\mu} \, y^j\bigg)\,
        d^\lam\ten d^\mu\ten \der_i
\\
    & \quad + h_{1}\,\del^i_j\, T\col \rho\rho\mu\,
        d^j\ten d^\mu\ten \der_i
\,,
\end{align*}
where all coefficients are real numbers.
\eLm

\bPf
We have the action of the group $W^{(1,1)}_{m,n} G$ on
the standard fiber $S_F=\Rn^{(m+n)*}\ten\Rn^{m*}\ten \Rn^n$
of $T^*\f E\ten T^*\f M\ten V\f E$
given by
\begin{align}
\bar \phi_\lam{}^i_\mu & = a^i_p\, \phi _\rho{}^p_\sig\,
        \tilde a^\rho_\lam\, \tilde a^\sig_\mu +
        a^i_p\, \phi _q{}^p_\sig\,
        \tilde a^q_{r\lam}\, a^r_s\, y^s\, \tilde a^\sig_\mu\,,
\quad
\bar \phi_j{}^i_\mu  = a^i_p\, \phi _q{}^p_\sig\,
        \tilde a^q_j\, \tilde a^\sig_\mu \,.
\end{align}

The corresponding equivariant mapping
$
\phi:\Rn^n\car \Rn^n\ten\Rn^{m*}\car T^s_m Q\car T^r_m R\to S_F
$
is given, by Corollary \ref{Cr2.3}, by
$$
\phi_A{}^i_\lam(y^i,y^i_\lam,\Lam \col \mu\lam{\nu,\f{\alp}},
        K\col ji{\mu,\f{\bet}}) =
        \psi_A{}^i_\lam(y^i,y^i_\lam,
        T\col \mu\lam{\nu;\f{\gam}},
        \widetilde R\col \mu\lam{\nu\kappa;\f{\del}},
        K\col ji\mu,R\col ji{\lam\mu;\f{\eps}})\,,
$$
where the multiindices $\f{\alp}, \dots, \f{\eps}$ are as in the proof
of Lemma \ref{Lm2.4}.

The equivariance with respect to fiber homotheties $(c\,\del^i_j)$
implies
$$
c\,\psi_\lam{}^i_\mu = \psi_\lam{}^i_\mu(c\, y^i,c\,y^i_\lam,
                T\col \mu\lam{\nu;\f{\gam}},
        \widetilde R\col \mu\lam{\nu\kappa;\f{\del}},
        K\col ji\mu,R\col ji{\lam\mu;\f{\eps}})\,.
$$
which gives, by the homogeneous function theorem, that
$\psi_\lam{}^i_\mu$ is linear in $y^i$ and $y^i_\lam$, i.e.
$$
\psi_\lam{}^i_\mu = f_{\lam\mu j}^i\, y^j +
        g_{\lam\mu j}^{i\nu}\, y^j_\nu\,.
$$

First, we shall discuss $g^{i\nu}_{\lam\mu j}$. We get, from the
equivariance with respect to base homotheties
$(c\, \del^\lam_\mu)$,
$$
c^{-1}\,g^{i\nu}_{\lam\mu j} = g^{i\nu}_{\lam\mu j}(
        c^{-(\parallel\f{\gam}\parallel+1)}\,
        T\col \mu\lam{\nu;\f{\gam}},
        c^{-(\parallel\f{\del}\parallel+2)}\,\widetilde
        R\col \mu\lam{\nu_1\nu_2;\f{\del}},
        c^{-1}\,K\col ji\lam,
        c^{-(\parallel\f{\eps}\parallel+2)}\,
        R\col ji{\lam\mu;\f{\eps}})\,
$$
which implies that $g^{i\nu}_{\lam\mu j}$ is polynomial with exponents
satisfying (\ref{Eq2.11}).
%$$
%-1 = - \sum^{s}_{i=0} (i+1)\,a_i - \sum^{s-1}_{j=0} (j+2)\, b_j
%        - c -
%         \sum^{r-1}_{k=0} (k+2)\, d_k \,.
%$$
We have the following 2 possible solutions:
$a_0=1$ and the other exponents vanish;
$c=1$ and the other exponents vanish.
Then
\begin{align*}
g^{i\nu}_{\lam\mu j} & =
  A^{i\nu\sig\tau}_{\lam\mu j\rho}\,
  T\col {\sig}{\rho}{\tau}
  + L^{i\nu q\sig}_{\lam\mu jp}\,
  K\col {q}{p}{\sig}
        = a_1\,\del^i_j\, T\col \lam\nu{\mu} + a_2\,\del^i_j\,
        \del_\lam^\nu\, T\col \rho\rho{\mu} + a_3\,\del^i_j\,
        \del_\mu^\nu\, T\col \lam\rho{\rho}
\\
&\quad +
        l_1\,\del^\nu_\lam\, K\col ji\mu
  +      l_2\,\del^\nu_\mu\, K\col ji\lam +
        l_3\,\del^i_j\,\del^\nu_\lam\, K\col pp\mu +
        l_4\,\del^i_j\,\del^\nu_\mu\, K\col pp\lam
\end{align*}

Similarly, for $f^{i}_{\lam\mu j}$ we get,
from the equivariance with respect to base homotheties
$(c\, \del^\lam_\mu)$,
$$
c^{-2}\,f^i_{\lam\mu j} = f^i_{\lam\mu j}(
                c^{-(\parallel\f{\gam}\parallel+1)}\,
        T\col \mu\lam{\nu;\f{\gam}},
        c^{-(\parallel\f{\del}\parallel+2)}\,\widetilde
        R\col \mu\lam{\nu_1\nu_2;\f{\del}},
        c^{-1}\,K\col ji\lam,
        c^{-(\parallel\f{\eps}\parallel+2)}\,
        R\col ji{\lam\mu;\f{\eps}})\,
$$
which implies that $f^i_{\lam\mu j}$ is polynomial with exponents
satisfying (\ref{Eq2.11}) with $-2$ on the left hand side.
%$$
%-2 = - \sum^{s}_{i=0} (i+1)\,a_i - \sum^{s-1}_{j=0} (j+2)\, b_j
%        - c -
%         \sum^{r-1}_{k=0} (k+2)\, d_k \,.
%$$
We have the following 6 possible solutions:
$a_0=2$ and the other exponents vanish;
$a_1=1$ and the other exponents vanish;
$a_0=1$, $c=1$ and the other exponents vanish;
$b_0=1$ and the other exponents vanish;
$c=2$ and the other exponents vanish;
$d_0=1$ and the other exponents vanish.
Then
\begin{align*}
f^i_{\lam\mu j} & =
  B^{i\sig_1\tau_1\sig_2\tau_2}_{\lam\mu j\rho_1\rho_2}\,
  T\col {\sig_1}{\rho_1}{\tau_1}\,T\col {\sig_2}{\rho_2}{\tau_2}
  + C^{i\sig\tau_1\tau_2}_{\lam\mu j\rho}\,
  T\col {\sig}{\rho}{\tau_1;\tau_2} +
  D^{i\sig\tau_1\tau_2}_{\lam\mu j\rho}\,
  R\col {\sig}{\rho}{\tau_1\tau_2}
\\
  & \quad   + P^{iq_1\tau_1 q_2\tau_2}_{\lam\mu jp_1p_2}\,
  K\col {q_1}{p_1}{\tau_1}\, K\col {q_2}{p_2}{\tau_2} +
  E^{iq\tau_1\tau_2}_{\lam\mu jp}\,
  R\col {q}{p}{\tau_1\tau_2} +
  N^{i\sig\tau_1 q\tau_2}_{\lam\mu j\rho p}\,
  T\col {\sig}{\rho}{\tau_1}\, K\col {q}{p}{\tau_2}
%\end{align*}
%where all coefficients are absolute invariant tensors, i.e.,
%\begin{align*}
%f^i_{\lam\mu j}
\\
& = b_1\,\del^i_j\, T\col \rho\rho{\lam}\,
        T\col \sig\sig{\mu} + b_2\,\del^i_j\, T\col \sig\rho{\lam}\,,
        T\col \rho\sig{\mu}+ b_3\,\del^i_j\, T\col \rho\rho{\sig}\,
        T\col \lam\sig{\mu} + c_1\,\del^i_j\, T\col \rho\rho{\lam;\mu}
\\
&\quad    + c_2\,\del^i_j\, T\col \rho\rho{\mu;\lam}
        + a_9\,\del^i_j\, T\col \lam\rho{\mu;\rho}
        + d_1\,\del^i_j\, R\col \rho\rho{\lam\mu}
        + d_2\,\del^i_j\, R\col \lam\rho{\rho\mu}
\\
&\quad         + p_1\,\del^i_j\, K\col pp\lam\, K\col qq\mu
        + p_2\,\del^i_j\, K\col qp\lam\, K\col pq\mu
        + p_3\, K\col ji\lam\, K\col pp\mu
\\
 &\quad   + p_4\, K\col ji\mu\, K\col pp\lam
        + p_5\, K\col pi\lam\, K\col jp\mu
        + p_6
\, K\col pi\mu\, K\col jp\lam
        + e_{1}\, R\col ji{\lam\mu}
        + e_{2}\,\del^i_j\, R\col pp{\lam\mu}
\\
&\quad         + n_1\, T\col \lam\rho\rho\, K\col ji\mu
        + n_2\, T\col \rho\rho\mu\, K\col ji\lam
        + n_{3}\, T\col \lam\rho\mu\, K\col ji\rho
        + n_{4}\,\del^i_j\, T\col \rho\rho\lam\, K\col pp\mu
\\
&\quad        + n_{5}\,\del^i_j\, T\col \rho\rho\mu\, K\col pp\lam
        + n_{6}\,\del^i_j\, T\col \lam\rho\lam\, K\col pp\rho\,.
\end{align*}

Further, from the equivariancy with respect to fiber homotheties,
we have
$$
\psi_j{}^i_\mu = \psi_j{}^i_\mu(c\, y^i,c\,y^i_\lam,
                T\col \mu\lam{\nu;\f{\gam}},
        \widetilde R\col \mu\lam{\nu\kappa;\f{\del}},
        K\col ji\mu,R\col ji{\lam\mu;\f{\eps}})\,.
$$
which implies, by the homogeneous function theorem, that
$\psi_j{}^i_\mu$ is independent of $y^i$ and $y^i_\lam$.

If we suppose the equivariance with respect to base homotheties
$(c\, \del^\lam_\mu)$. We have
$$
c^{-1}\,\psi_j{}^i_\mu = \psi_j{}^i_\mu(
                c^{-(\parallel\f{\gam}\parallel+1)}\,
        T\col \mu\lam{\nu;\f{\gam}},
        c^{-(\parallel\f{\del}\parallel+2)}\,\widetilde
        R\col \mu\lam{\nu_1\nu_2;\f{\del}},
        c^{-1}\,K\col ji\lam,
        c^{-(\parallel\f{\eps}\parallel+2)}\,
        R\col ji{\lam\mu;\f{\eps}})\,
$$
which implies that $\psi_j{}^i_\mu$ is polynomial with exponents
satisfying  the equation (\ref{Eq2.11}),
%$$
%-1 = - \sum^{s}_{i=0} (i+1)\,a_i - \sum^{s-1}_{j=0} (j+2)\, b_j
%        - c -
%         \sum^{r-1}_{k=0} (k+2)\, d_k \,,
%$$
i.e.
$$
\psi_j{}^i_\mu  = H^{i\sig\tau}_{j\mu\rho}\, T\col \sig\rho\tau
        + M^{iq\rho}_{j\mu p}\, K\col qp\rho
        = h_1\,\del^i_j\, T\col \rho\rho\mu +
        m_{1}\, K\col ji\mu + m_{2}\,\del^i_j\, K\col pp\mu\,.
$$
%where $H^{i\sig\tau}_{j\mu\rho}, M^{iq\rho}_{j\mu p}$ are absolute
%invariant tensors. Then
%$$
%\psi_j{}^i_\mu =
%$$

Finally, the equivariance of $\psi_j{}^i_\mu$ with respect to
elements of the type $(\del^i_j, \del ^\lam_\mu, a^i_{j\lam})$
implies
$m_{1}=m_{2}=0$ and we have
\bEq\label{Eq3.1}
\psi_j{}^i_\mu = h_1\,\del^i_j\, T\col \rho\rho\mu\,
\eEq
and the equivariance of $\psi_\lam{}^i_\mu$ with respect to
elements of the type $(\del^i_j, \del ^\lam_\mu, a^i_{j\lam})$
implies
$p_i=0$, $i=1,\dots,10$, $n_4=n_5=n_6=0$, $n_1=-a_3$, $n_2=-a_2-h_1$,
$n_3=-a_1$ and the other coefficients are arbitrary. Then
\begin{align}\label{Eq3.2}
\psi_\lam{}^i_\mu & = \phi_\lam{}^i_\mu
        = g^{i\nu}_{\lam\mu j}\, y^j_\nu + f^{i}_{\lam\mu j}\, y^j
\\
  & = a_1\, T\col \lam\nu\mu\, y^i_\nu
        +  a_2\, T\col \rho\rho\mu \, y^i_\lam
        + a_3\, T\col \lam\rho\rho \, y^i_\mu
        + b_1\, y^i\, T\col \rho\rho\lam \, T\col \sig\sig\mu
        + b_2\, y^i\, T\col \sig\rho\lam \, T\col \rho\sig\mu
        \nonumber
\\
&\quad         + b_3\, y^i\, T\col \rho\rho\sig \, T\col \lam\sig\mu
      + c_1\, y^i\, T\col \rho\rho{\lam;\mu}
        + c_2\, y^i\, T\col \rho\rho{\mu;\lam}
        + c_3\, y^i\, T\col \lam\rho{\mu;\rho}
        \nonumber
 \\
&\quad        - a_3\, T\col \lam\rho{\rho}\, K\col ji\mu \, y^j
        - (a_2+a_{10})\, T\col \rho\rho{\mu}\, K\col ji\lam \, y^j
        - a_1\, T\col \lam\rho{\mu}\, K\col ji\rho \, y^j
        \nonumber
\\
   &\quad + d_1 \, y^i\, R\col \rho\rho{\lam\mu}
        + d_2 \, y^i\, R\col \lam\rho{\rho\mu}
        + e_1 \,y^i\, R\col pp{\lam\mu}
        + e_2 \, R\col ji{\lam\mu} \, y^j\,.\nonumber
\end{align}

Summerizing (\ref{Eq3.1}) and (\ref{Eq3.2}) we get Lemma \ref{Lm3.6}.
\ePf

\bTh\label{Th3.9}
All natural operators transforming a linear connection $K$ on $\f E$
and a classical linear connection $\Lam$ on $\f M$ into
connections on $\pi^1_0:J_1\f E\to \f E$ are of the maximal order one
and form the 14-parameter family
$$
\widetilde{\Gam}(\Lam,K) = \chi(D(\Lam,K)) + \phi(\Lam,K)\,,
$$
where $\chi(D(\Lam,K))$ is the connection given
by Theorem \ref{Th3.4}
and $\phi(\Lam,K)$ is the 14-parameter family of natural tensor
fields given by Lemma \ref{Lm3.6}.
\hfill\qedsymbol
\eTh

Let us recall that we have
the natural complementary contact maps
\beq
{\cyrm d}:J^1\f E\ucar{\f E} T\f M\to T\f E \,,
\quad
\theta:J^1\f E\ucar{\f E}T\f E\to V\f E \,,
\eeq
with the coordinate expressions
\bEq\label{de}
{\cyrm d} =  d^\lam \ten (\der_\lam+y^i_\lam\der_i) \,,
\quad
\theta = (d^i-y^i_\lam d^\lam) \ten \der_i \,.
\eEq
Then we have the following, by using the notation of Section 2,
geometric description of $\widetilde\Gam(\Lam, K)$.

\bTh\label{Th3.7}
All connections on $J^1\f E$ naturally given by
$\Lam$ (in order $s$) and by
$K$ (in order $r$,
$s\ge r-2$) are of the maximal
order one  and are of the  form
\begin{align*}
\widetilde{\Gam}(\Lam,K)=\Gam(\Lam,K) +
        \theta\com h^K\big(S(\Lam)\big)
        +L\ten G(\Lam,K) + e_2\, R[K](L)
       + h_1\, \nu_K \ten \hat T\,.
\end{align*}
\vglue-1.3\baselineskip{\ }\hfill\END
\eTh

\bRm\label{Rm3.8}
The 14-parameter family $\widetilde\Gam(\Lam, K)$ from Theorem
\ref{Th3.9} can be obtained
from the 15-parameter family of Theorem \ref{Th2.7} by applying the
operator $\chi$, i.e.
$\widetilde\Gam(\Lam, K)=\chi\big(\widetilde D(\Lam,K)\big)$.
In fact $\chi\big(\widetilde D(\Lam,K)\big)=
\chi\big(D(\Lam,K)\big) +\widetilde{\chi}\big(\Phi(\Lam,K)\big)$,
where
$$
  \widetilde{\chi}\equiv\id_{T^*\f E}
        \ten \theta \ten \K{d}:J^1\f E\ucar{\f E}
        T^*\f E\ten T\f E\ten T^*\f E\to
        T^*\f E\ten V\f E\ten T^*\f M\,.
$$
But
$$
\widetilde{\chi}
\big(h^K(\hat T\ten I_{T\f M})\big)
=- \widetilde{\chi}
(\hat T\ten \nu_{K})
$$
which implies that operators standing in the family of Lemma \ref{Lm2.4}
with coefficients $a_3$ and $h_{2}$ admit the same operator,
standing with the coefficient $a_3$ in the family
of Lemma \ref{Lm3.6}.
It is the reason why $\widetilde\Gam(\Lam,K)$ is only 14-parameter family.
\hfill\qedsymbol
\eRm

\bRm\label{Rm3.10}
From the coordinate expression of $\Gam(\Lam,K)$ and $\phi(\Lam,K)$
it is easy to see that all connections $\widetilde\Gam(\Lam,K)$ are
affine.
\hfill\qedsymbol
\eRm

\bCr\label{Cr3.11}
All natural operators transforming a  general linear
connection $K$ on $\f E$
and a symmetric classical connection $\Lam$ on
$\f M$ into  connections on $J^1\f E$
are of the maximal order one
and form the 4-parameter family
$$
\widetilde{\Gam}(\Lam,K)
        =  \Gam(\Lam,K) + L\ten\big(d_1\, C^1_1 R[\widetilde\Lam]
        + d_2\,  C^1_2 R[\widetilde\Lam]+
         e_1\, C^1_1 R[K]\big) + e_2\, R[K](L)\,.\eqno{\qedsymbol}
$$
\eCr

%--------------------------------------------------------------------%

%--------------------------------------------------------------------%

\end{document}